\documentclass[11pt]{article}
\usepackage{amssymb}

\def\F{{\cal F}}

\def\M{{\cal M}}
\def\N{{\cal N}}

\newtheorem{theorem}{Theorem}
\newtheorem{prop}[theorem]{Proposition}
\newtheorem{lemma}[theorem]{Lemma}

\begin{document}

\bigskip
\centerline{{\Large\bf Frobenius submanifolds}}

\vspace{.3in}
\centerline{{\bf I.A.B. Strachan}}\vspace{.1in}
\centerline{Department of Mathematics, University of Hull,}
\vspace{.1in}
\centerline{Hull, HU6 7RX, England.}
\vspace{.1in}
\centerline{e-mail: i.a.b.strachan@hull.ac.uk}

\vspace{.4in}
\centerline{{\bf Abstract}}

\vspace{.3in}
\small
\hskip 10mm\parbox{5.0in}{The notion of a Frobenius submanifold -- a submanifold of a Frobenius
manifold which is itself a Frobenius manifold with respect to structures induced from
the original Frobenius manifold -- is studied. Two dimensional submanifolds are particularly
simple. More generally, sufficient conditions are given for a submanifold to be a so-called
natural Frobenius submanifold. These ideas are illustrated using examples of Frobenius manifolds
constructed from Coxeter groups, and for the Frobenius manifolds governing the quantum
cohomology of $\mathbb{CP}^2$ and $\mathbb{CP}^1 \times \mathbb{CP}^1\,.$}
\normalsize

\bigskip

\section{Introduction}

Substructures abound within mathematics. The purpose of this paper is to introduce the
notion of a Frobenius submanifold -- a submanifold of a Frobenius manifold which
is itself a Frobenius manifold with respect to structures induced from the original
Frobenius manifold. Certain specialized examples have appeared in the literature before, but
the approach was more algebraic than geometric, the submanifolds being hyperplanes [Z].
The paper is laid out as follows. In section 2 a more general framework of induced substructures
is given, with Frobenius submanifolds being introduced in section 3. So called natural
Frobenius submanifolds are studied in more detail in section 4, and in the remaining sections
a series of examples based on the foldings of Coxeter graphs and on the
quantum cohomology of certain projective spaces are studied.

\section{Submanifolds and their induced structures}

Let $\M$ be some manifold endowed with a metric $\eta=<,>\,.$ Suppose further that on each
tangent space $T_t\M$ one has a commutative multiplication of vectors
\[
\circ\,:\quad\quad T_t\M \times T_t\M \longrightarrow T_t\M\,,
\]
varying smoothly over the manifold. Moreover, it will be assumed that this multiplication
is compatible with the metric, in the sense that
\[
<a\circ b,c> = <a,b\circ c> \quad\quad \forall a\,,b\,,c\,\in T_t\M\,.
\]
This property is known as the Frobenius condition. Let $\F$ denote the triple
$\F=\{\M,\eta,\circ\}\,.$ This will be called a Frobenius structure.

Let $\N\subset\M$ be a submanifold of $\M\,.$ One may defined an induced $\F$ structure on
$\N\,,$ denoted by $\F_\N=\{\N,\eta_N,\star\}\,,$ as follows. The metric $\eta_\N=<,>_\N$
is just the induced
metric on $\N\,,$ and $\star$ is defined by
\[
a\star b = pr(a \circ b) \quad \forall a\,,b\,\in T_\tau\N \subset T_\tau\M\,,
\]
where $pr$ denotes the projection (using the original metric $\eta$ on $\M$) of
$a\circ b\in T_\tau\M$ onto $T_\tau\N\,.$ This induced multiplication may have very
different algebraic properties than those of its progenitor.

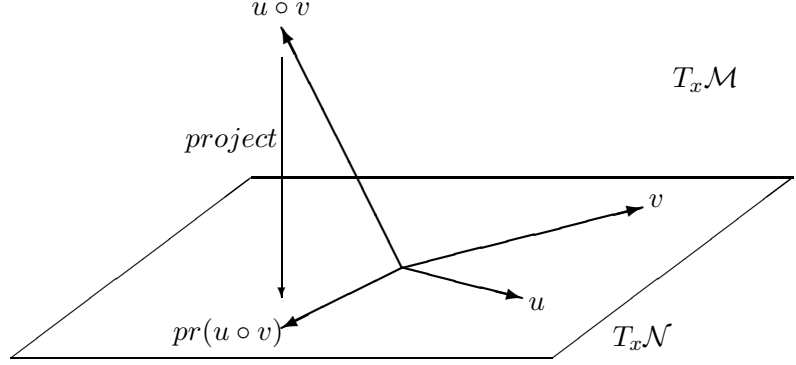
\begin{figure}

\setlength{\unitlength}{0.8cm}
\begin{picture}(17,7)
\put(3,0){\line(1,0){9}}
\put(12,0){\line(4,3){4}}
\put(3,0){\line(4,3){4}}
\put(7,3){\line(1,0){9}}
\thicklines
\put(9.5,1.5){\vector(4,-1){2}}
\put(9.5,1.5){\vector(4,1){4}}
\put(9.5,1.5){\vector(-1,2){2}}
\put(9.5,1.5){\vector(-2,-1){2}}
\put(11.6,0.8){$u$}
\put(13.6,2.5){$v$}
\put(7.0,5.7){$u\circ v$}
\thinlines
\put(7.5,5.0){\vector(0,-1){4}}
\put(5.7,0.3){$pr(u\circ v)$}
\put(5.9,3.5){$project$}
\put(13,0.2){$T_x{\cal N}$}
\put(14,4.5){$T_x{\cal M}$}
\end{picture}

\caption{The definition of the induced multiplication}
\end{figure}

\begin{lemma}
The induced structure $\F_\N$ satisfies the Frobenius condition
\[
<a\star b,c>_\N = <a,b\star c>_\N \quad\quad \forall a\,,b\,c\,\in T_\tau\N\,.
\]
Hence $\F_\N$ is a Frobenius structure.
\end{lemma}

\noindent The proof follows immediately from the definitions. An alternative proof will be
given below. Before this some general results will be given; this will also serve
to fix the notation that will subsequently be used in this paper.

Let $t^i\,,i=1\,,\ldots \, m=dim\M$ be local coordinates on $\M\,.$ With these the submanifold
$\N$ may be defined by the parametrization
\begin{equation}
t^i = t^i(\tau^\alpha)\,,\quad \alpha=1\,,\ldots \, n=dim\N\,,i=1\,,\ldots \, m=dim\M\,,
\label{basisN}
\end{equation}
and so a basis for $T_\tau\N$ is given by
\[
\frac{\partial ~}{\partial\tau^\alpha} = \frac{\partial t^i}{\partial \tau^\alpha}
\frac{\partial ~}{\partial t^i}\,.
\]
In these coordinates the induced metric on $\N$ is given by\footnote{The symbol $\eta$ will
be used to denote a metric on either $\M$ or $\N\,,$ with Greek indices denoting structures
on $\N$ and Latin indices structures on $\M\,.$ This convention will be used throughout this
paper.}
\begin{equation}
\eta_{\alpha\beta}=
\frac{\partial t^i}{\partial \tau^\alpha}
\frac{\partial t^j}{\partial \tau^\beta}\eta_{ij}\,,
\label{inducedmetric}
\end{equation}
where $\eta_{ij}$ is the metric on $\M\,.$ The basis (\ref{basisN}) may be extended to a basis
for $T_t\M\,,$ so
\begin{equation}
\frac{\partial~}{\partial t^i} = A^\alpha_i \frac{\partial~}{\partial \tau^\alpha} +
n_i^{\tilde\alpha} \frac{\partial~}{\partial\nu^{\tilde\alpha}}   \,,
\label{decomposition}
\end{equation}
where ${\tilde\alpha}=1\,,\ldots\,,m-n\,$ and
\[
\frac{\partial~}{\partial\nu^{\tilde\alpha}} \in (T_\tau\N)^\perp\,.
\]
Using the metrics on $T_t\M$ and $T_n\N$ one
obtains
\[
A^\alpha_i = \eta^{\alpha\beta}\eta_{ij} \frac{\partial t^j}{\partial \tau^\beta}\,.
\]
The multiplication on $T_t\M$ may be defined in terms of a set of structure functions
$c_{ij}^{~~~k}(t^r)\,:$
\[
\frac{\partial~}{\partial t^i} \circ \frac{\partial~}{\partial t^j}=
c_{ij}^{~~~k} \frac{\partial~}{\partial t^k}\,.
\]
With these one may find the induced structure functions for the multiplication on $T_\tau\N\,.$
\begin{eqnarray*}
\frac{\partial~}{\partial\tau^\alpha} \circ \frac{\partial~}{\partial\tau^\beta}
& = & \frac{\partial t^i}{\partial \tau^\alpha}\frac{\partial t^j}{\partial \tau^\beta}
\left.c_{ij}^{~~~k}
\right|_\N \frac{\partial~}{\partial t^k}\,, \\
& = & \frac{\partial t^i}{\partial \tau^\alpha}\frac{\partial t^j}{\partial \tau^\beta}
\left.c_{ij}^{~~~k}
\right
|_\N \Big[
A^\gamma_k \frac{\partial~}{\partial \tau^\gamma} +
n_k^{\tilde\gamma} \frac{\partial~}{\partial\nu^{\tilde\gamma}}\Big]\,.
\end{eqnarray*}
Hence
\begin{eqnarray*}
\frac{\partial~}{\partial\tau^\alpha} \star \frac{\partial~}{\partial\tau^\beta}
& = & pr \Big[
\frac{\partial~}{\partial\tau^\alpha} \circ \frac{\partial~}{\partial\tau^\beta}
\Big]\,,\\
& = & \frac{\partial t^i}{\partial \tau^\alpha}\frac{\partial t^j}{\partial \tau^\beta}
\left.c_{ij}^{~~~k}\right|_\N 
A^\gamma_k \frac{\partial~}{\partial \tau^\gamma}\,,\\
& = & c_{\alpha\beta}^{~~~\gamma} \frac{\partial~}{\partial\tau^\gamma}\,,
\end{eqnarray*}
where the induced structure functions are given by
\begin{equation}
c_{\alpha\beta}^{~~~\gamma}=
\frac{\partial t^i}{\partial\tau^\alpha}
\frac{\partial t^j}{\partial\tau^\beta}
\frac{\partial t^r}{\partial\tau^\delta}
\eta_{kr}\eta^{\gamma\delta}\left.c_{ij}^{~~~k}\right|_\N\,.
\label{inducedmult}
\end{equation}

\bigskip

\noindent{\bf Proof~} With the notion set up the proof of the proposition is straightforward.
The Frobenius property on $\M$ is equivalent to the condition that the tensor
\[
c_{ijk}=\eta_{kl}c_{ij}^{~~~l}
\]
is totally symmetric (recall that $\circ$ is, by definition, commutative). It follows
from this and (\ref{inducedmult}) that
\begin{equation}
c_{\alpha\beta\gamma} =
\frac{\partial t^i}{\partial\tau^\alpha}
\frac{\partial t^j}{\partial\tau^\beta}
\frac{\partial t^k}{\partial\tau^\gamma}
\left.c_{ijk}\right|_\N
\label{inducedthirdderv}
\end{equation}
is also totally symmetric. Hence the induced structure $\F_\N$ inherits the Frobenius
property.


\bigskip

\noindent{\bf Example} Consider the Jordan algebra defined by the
\begin{eqnarray*}
e_1 \circ e_i & = & + e_i \,, \quad i=1\,,\ldots \,, m\,,\\
e_i \circ e_i & = & - e_1 \,, \quad i=2\,,\ldots \,, m\,,\\
e_i \circ e_j & = & 0 \quad\quad\quad {\rm otherwise\,.}
\end{eqnarray*}
One may show that with the inner product defined by $\eta_{ij}=c_{ij}^{~~~k}c_{kl}^{~~~l}$
(where $c_{ij}^{~~~k}$ are the structure constants of this algebra) this algebra
has the Frobenius property [S1]. These may then be used to define a trivial $\F$-structure -
trivial in the sense that the structures do not vary are the tangent space varies. The
above proposition may then be used to find examples of other, non-trivial, $\F$-structures.

\bigskip

\noindent In what follows the idea of a natural substructure will be important.

\bigskip

\noindent{\bf Definition} A substructure $\F_\N$ of a Frobenius structure $\F$ is said to
be {\it natural} if
\[
a\star b = a \circ b \,, \quad \forall\,a\,,b\,\in T_\tau\N\,,
\]
that is, no projection onto $T_\tau\N$ is required, for all points $x\in \N\,.$

\bigskip

\noindent In terms of the local coordinates, this means that the $n(n+1)(m-n)/2$ conditions
$\Xi_{\alpha\beta}^{~~~\tilde\gamma}$ must vanish, where
\begin{equation}
\Xi_{\alpha\beta}^{~~~\tilde\gamma}=
\frac{\partial t^i}{\partial\tau^\alpha}
\frac{\partial t^j}{\partial\tau^\beta}
\left.c_{ij}^{~~~k}\right|_\N n_k^{~{\tilde\gamma}}\,.
\label{obstruction}
\end{equation}

\bigskip

\noindent{\bf Example} Let $I\subset\{1\,,2\,,\ldots\,,m\}\,$ and suppose that $\N$ is given by the conditions
$t^i=0$ for $i\notin I\,.$ Then the obstruction reduces to the algebraic condition
\[
\left.c_{ij}^{~~~k}\right|_N = 0 \,,\quad\quad i\,,j \in I\,, k\notin I\,.
\]
This condition was derived in [Z] in the context of Frobenius manifolds constructed
from Coxeter groups (see section 5). Here it is a
specialization of the more general condition (\ref{obstruction}).

\bigskip

\section{Frobenius manifolds}

One particular class of Frobenius structures are Frobenius manifolds.
A Frobenius manifold may be defined as follows [D]. Let $F=F(t^i)$ be a
function -- the prepotential --
defined on some region
$\M\subset {\mathbb{R}}^m$ (sometimes $\M\subset {\mathbb{C}}^m$)
such that the third derivatives
\[
c_{ijk}=\frac{\partial^3 F}{\partial t^i \partial t^j \partial t^k}
\]
satisfy the following conditions:

\begin{itemize}

\item[$\bullet$] Normalization:

\[
\eta_{ij}=c_{1ij}
\]
is a constant, nondegenerate matrix. Let $\eta^{ij}=(\eta_{ij})^{-1}\,.$ These may be used
to raise and lower indices.

\item[$\bullet$] Associativity: the functions

\[
c_{ij}^{~~~k} = \eta^{kl} c_{ijl}
\]
define an associative, commutative algebra
\[
\frac{\partial~}{\partial t^i}\circ
\frac{\partial~}{\partial t^j} = c_{ij}^{~~~k}
\frac{\partial~}{\partial t^k}
\]
on each tangent space $T_t\M$ with unity element $\mathbb{I}\,,$ so
$\mathbb{I}\circ a = a \quad\forall a\in T_t\M\,.$ The above normalization implies
that $\mathbb{I}=\partial_{t^1}\,.$ The resulting differential equation for the
prepotential is known as the Witten-Dijkgraaf-Verlinde-Verlinde (or WDVV) equation.

\item[$\bullet$] Homogeneity: The function $F$ must be quasi-homogeneous, so
\[
\mathcal{L}_E = d_F \, F + {\rm quadratic~terms}
\]
where $\mathcal{L}_E$ is the Lie derivative along the Euler vector field
\[
E=(q^i_j t^j + r^i) \frac{\partial ~}{\partial t^i}\,.
\]

\end{itemize}
The most common form (which is canonical under certain additional requirements)
for $\eta_{ij}$ is the antidiagonal form
\[
\eta_{ij} = \delta_{i+j,m+1}\,,
\]
and this form will be assumed throughout this paper. It then follows from the
above axioms that the prepotential takes the general form\footnote{To avoid a plethora
of brackets in terms such as $(t^2)^2(t^3)^3$ the indices on $t$ will be written downstairs
in {\sl explicit} formulae, so $t_i=t^i\,,$ not $t_i=\eta_{ij}t^j\,.$}
\begin{equation}
F = \frac{1}{2} t_1^2 t_m + \frac{1}{2} t_1 \sum_{j=2}^{m-1} t_i t_{m-i+1}
+ f(t_2\,,\ldots\,,t_m)\,.
\label{canonicalF}
\end{equation}

\noindent It will be assumed that the Euler vector field $E$ takes the form
\[
E=\sum_i d_i t^i \frac{\partial ~}{\partial t^i} +
\sum_{i|d_i=0} r^i \frac{\partial ~}{\partial t^i}
\]
with $d_1=1\,,$ and with the canonical form (\ref{canonicalF}) for the prepotential
\[
q_i+q_{m+1-i}=d\,,
\]
where $d=3-d_F$ and $d_i=1-q_i\,.$

\bigskip

\noindent{\bf Example} $m=2\,.$ The equations of associativity are vacuous, so any function
\[
F=\frac{1}{2} t_1 t_2^2 + f(t_2)
\]
defines a Frobenius manifold. If the quasihomogeneity condition is now used the otherwise
free function $f(t_2)$ is constrained to take one of five forms.

\bigskip

\noindent{\bf Example} $m=3\,.$ The equations of associativity results in a single
differential equation for $f(x,y)\,,$
\[
f_{xxy}^2 = f_{yyy} + f_{xxx} f_{xyy}\,.
\]
If the quasihomogeneity condition is now used this equation may be reduced to various
third order
ordinary differential equation, each equivalent to a Painlev\'e VI equation.

\bigskip

On a submanifold $\N$ one may, as well as the induced $\F_\N$ structures,
define an induced vector field
\[
E_\N = pr \left. E\right|_\N\,.
\]
This raises a numer of questions on whether the induced structures are quasi-homogeneous
with respect to the induced Euler vector field, and in particular:

\begin{itemize}

\item[$\bullet$] For what families of submanifolds does
\[
E_\N = (q^\alpha_\beta \tau^\beta + r^\alpha)\frac{\partial ~}{\partial \tau^\alpha}\,,
\]
since in general $E_\N$ will not be linear in $\tau^i\,?$

\item[$\bullet$] For what families of submanifolds does
\[
E_\N = \left.E\right|_\N\,
\]
or equivalently, $(\left.E\right|_\N)^\perp=0\,$?

\end{itemize}
It will be shown in section 3 that for natural Frobenius submanifolds, the
second condition implies the first.

\bigskip

\bigskip

\noindent{\bf Definition} Let $\F$ be a Frobenius manifold. A submanifold $\N$ be said
to be a Frobenius submanifold if $\F_\N$ is a Frobenius manifold with respect to the
induced structures.

\bigskip

\bigskip

\noindent{\bf Definition} A natural Frobenius submanifold $\N$ is a Frobenius submanifold
where
\[
a \star b = a \circ b \,,\quad \forall a\,,b\,\in T_\tau\N\,,
\]
or equivalently, $\big(  ( a\star b )|_\N \big)^\perp=0\,.$

\bigskip

\noindent For the rest of this section the quasihomogeneity condition will be ignored,
concentrating instead on properties of the induced multiplication on
two dimensional submanifolds. It turns out that
two dimensional Frobenius submanifolds are particularly simple.

\medskip

\begin{prop} Let $\F=\{\M,\eta,\circ\}$ be a Frobenius manifold and let $\N$ be
a two dimensional submanifold such the unity vector field at all points of $\N$
is always tangent to $\N\,.$ Then $\F_\N$ is a Frobenius manifold.
\end{prop}

\bigskip

\noindent{\bf Proof~} The proof will only be given for $dim\,\M=3\,,$ the general case being
a direct generalization of the lower dimensional result. To fulfil the tangential condition
the surface may be parametrized by
\[
\pmatrix{t_1 \cr t_2 \cr t_3} = \pmatrix{1\cr 0 \cr 0} \tau_1 +
\pmatrix{a(\tau_2) \cr b(\tau_2) \cr c(\tau_2) }\,,
\]
this ensuring that $\partial_{t_1}=\partial_{\tau_1}\,.$ The induced metric on the ruled
surface is automatically flat, and flat coordinates are given by
\[
\pmatrix{t_1 \cr t_2 \cr t_3} = \pmatrix{1\cr 0 \cr 0} \tau_1 +
\pmatrix{-\frac{1}{2} \int b_{\tau_2}^2 \, d\tau_2 \cr b(\tau_2) \cr \tau_2}
\]
in which the induced metric is just $\left.\eta\right|_N = 2 d \tau_1 \, d\tau_2\,.$

In this parametrization
\begin{eqnarray*}
\frac{\partial ~}{\partial \tau_1} & = & +\frac{\partial ~}{\partial t_1}\,, \\
\frac{\partial ~}{\partial \tau_2} & = & -\frac{1}{2} \frac{\partial ~}{\partial t_1}
+ b_{\tau_2} \frac{\partial ~}{\partial t_2} + \frac{\partial ~}{\partial t_3}\,.
\end{eqnarray*}
Since $\partial_{t_1}$ is the unity vector field
\begin{eqnarray*}
\frac{\partial ~}{\partial \tau_1} \star \frac{\partial ~}{\partial \tau_1}
& = & \frac{\partial ~}{\partial \tau_1} \,, \\
\frac{\partial ~}{\partial \tau_1} \star \frac{\partial ~}{\partial \tau_2}
& = & \frac{\partial ~}{\partial \tau_1} \, \\
\end{eqnarray*}
and so it just remains to calulate $\partial_{\tau_2} \circ \partial_{\tau_2}$ and project
this onto $T_x \N\,.$ The vector
\[
\frac{\partial ~}{\partial\nu} = \frac{\partial ~}{\partial t_2} -
b_{\tau_2} \frac{\partial ~}{\partial \tau_1}
\]
is perpendicular to $T_\tau\N$ (it is not necessary to normalise here) and hence
\begin{eqnarray*}
\frac{\partial ~}{\partial \tau_2} \circ \frac{\partial ~}{\partial \tau_2}
& = &
\phantom{+}\Bigg(- \frac{3}{4} b_{\tau_2}^4 + 
b_{\tau_2}^3 \left.c_{222}\right|_\N + 3
b_{\tau_2}^2 \left.c_{223}\right|_\N + 3
b_{\tau_2} \left.c_{233}\right|_\N  + \left.c_{333}\right|_\N \Bigg)
\frac{\partial ~}{\partial \tau_1}   \,, \\
& &
+\Bigg( -b_{\tau_2}^3 + 
b_{\tau_2}^2 \left.c_{222}\right|_\N + 2
b_{\tau_2} \left.c_{223}\right|_\N + 
\left.c_{233}\right|_\N \Bigg)  
\frac{\partial ~}{\partial \nu}\,.
\end{eqnarray*}
Hence, on projecting onto $T_\tau\N\,,$
\begin{eqnarray*}
\frac{\partial ~}{\partial \tau_2} \star \frac{\partial ~}{\partial \tau_2}
& = & pr\Bigg(
\frac{\partial ~}{\partial \tau_2} \circ \frac{\partial ~}{\partial \tau_2}
\Bigg) \,, \\
& = & [{\rm ~function~of~}\tau_2 \,]\, \frac{\partial ~}{\partial \tau_1}
\,.
\end{eqnarray*}
If $dim\,\N$ was greater than two one would now have to check that this multiplication
was associative, but in two dimensions the associativity condition is vacuous, and this
induced structure is automatically a Frobenius submanifold with prepotential
\[
F_N = \frac{1}{2} \tau_1^2\tau_2 + \int\!\!\!\int\!\!\!\int
[{\rm ~function~of~}\tau_2] \, d \tau_2\,d \tau_2 d \tau_2\\.
\]
The condition the surface to be a natural Frobenius submanifold is thus
\[
b_{\tau_2}^3 = 
b_{\tau_2}^2 \left.c_{222}\right|_\N + 2
b_{\tau_2} \left.c_{223}\right|_\N + 
\left.c_{233}\right|_\N\,,
\]
a first order ordinary differential equation of degree three. Note that in general
\[
\left.F\right|_N \neq F_\N\,.
\]
Thus any two dimensional manifold ruled in this way is a Frobenius submanifold.

\bigskip

\bigskip

\section{Natural Frobenius submanifolds}

In this section sufficient conditions will be derived to ensure that a flat submanifold
of a Frobenius manifold is a natural Frobenius submanifold.

\begin{theorem} Let $\N$ be a flat submanifold of a Frobenius manifold $\M$ with
\begin{eqnarray*}
\Big( {\mathbb{I}}|_\N \Big)^\perp & = & 0 \,, \\
\Big(  ( a \circ b)|_\N \Big)^\perp & = & 0 \,, \quad \forall a\,,b\in T_\tau\N\,,\\
\Big( E|_\N \Big)^\perp & = & 0 \,.
\end{eqnarray*}
Then $\N$ is a natural Frobenius submanifold
 
\end{theorem}

\medskip

\noindent With so many conditions on $\N$ the result might seem inevitable, but it is not clear
that a prepotential exists, or that the induced Euler vector field is linear, or that the
induced prepotential is quasihomogeneous with repect to the induced Euler vector field.

\bigskip

\noindent{\bf Proof~} Since $\N$ is flat one may, by solving the Gauss-Manin equations,
find coordinates so that the components of the induced metric (\ref{inducedmetric}) are
constant - the so-called flat coordinates. The geometric properties of a
flat submanifold in a flat
manifold is summarized in the appendix.

\medskip

\noindent{\underline {Existence of induced prepotential}}

\medskip

Since $a\circ b=a \star b$ it follows immediately that
$\circ$ is a commutative, associative multiplication with induced structure
functions given by (\ref{inducedmult})\,. 
The existence of an induced prepotential $F_\N$ such that
\[
c_{\alpha\beta\gamma}=
\frac{\partial^3 F_\N}{\partial \tau^\alpha \partial \tau^\beta \partial \tau^\gamma}\,
\]
is given by the integrability conditions
\[
\frac{\partial c_{\alpha\mu\nu}}{\partial\tau^\beta} -
\frac{\partial c_{\beta\mu\nu}}{\partial\tau^\alpha} = 0\,.
\]
Using (\ref{inducedthirdderv}) and (\ref{a4})
\begin{eqnarray*}
\frac{\partial c_{\alpha\mu\nu}}{\partial\tau^\beta} -
\frac{\partial c_{\beta\mu\nu}}{\partial\tau^\alpha} & = &
\sum_{\rm similar~terms} \pm
\frac{\partial t^i}{\partial\tau^\sigma}
\frac{\partial^2 t^j}{\partial\tau^\beta\partial\tau^\mu}
\frac{\partial t^k}{\partial\tau^\nu}
\left.c_{ijk}\right|_\N\,, \\
& = & \sum_{\rm similar~terms} \pm
\frac{\partial t^i}{\partial\tau^\sigma}
\frac{\partial t^k}{\partial\tau^\nu}
\Omega^{~~~\tilde\alpha}_{\beta\mu} n_{\tilde\alpha}^{~ j} \left.c_{ijk}\right|_\N\,, \\
& = & \sum_{\rm similar~terms} \pm\Omega^{~~~\tilde\alpha}_{\beta\mu} \Xi_{{\tilde\alpha}\alpha\nu}\,.
\end{eqnarray*}
The obstruction to the existence of a prepotential is thus
\[
{\rm obstruction} = \Omega^{[\alpha}_{(\mu} \Xi^{\beta]}_{\nu)}
\]
(supressing the sum over $\tilde\alpha$). Two simple cases where this obstruction
vanishes are:
\begin{eqnarray*}
\Xi & = & 0 \,,\\
\Xi & = & \Omega\,.
\end{eqnarray*}
Thus on a natural submanifold these obstructions vanish (since the $\Xi$ vanish)
and an induced prepotential
$F_\N$ exists. Proposition [2] shows that this condition is not necessary.

\medskip

\noindent{\underline {Existence of unity element}}

\medskip

Recall that
\[
\mathbb{I}=\frac{\partial~}{\partial t^1}\,.
\]
Using  this and (\ref{decomposition}) it follows that $n_1^{~{\tilde\alpha}}=0\,,$ and this
together with (\ref{a4}) implies that
\[
\frac{\partial^2 t^m}{\partial\tau^\alpha\partial\tau^\beta}=0\,,\quad\quad(m={\rm dim\,}\M)\,.
\]
Hence $t^m=\mu_\alpha \tau^\alpha + \beta\,,$ where $\mu_\alpha$ and $\beta$ are constants.
Linear transformations, which do not affect the flatness of the $\tau$-coordinates, may
be used to fix $t^m=\tau^n\,.$ This ensures that
\begin{eqnarray*}
\mathbb{I}_\N & = & pr \big( \mathbb{I}|_\N \big )\,,\\
& = & \frac{\partial t^m}{\partial \tau^\alpha} \eta^{\alpha\beta}
\frac{\partial ~}{\partial\tau^\beta}\,,\\
& = & \frac{\partial ~}{\partial\tau^1}\,.
\end{eqnarray*}
The parametrization of the submanifold must have the generic form
\begin{equation}
\begin{array}{rcl}
t^1 & = & \tau^1 + f^1(\tau^2\,,\ldots\,,\tau^n)\,, \\
t^i & = & f^i(\tau^2\,,\ldots\,,\tau^n)\,\quad i=2\,,\ldots\,,m-1\,,\\
t^m & = & \tau^n\,.
\end{array}
\label{param}
\end{equation}
Having set up appropriate coordinates on $\N$ the required normalization on the submanifold
is straightforward:
\begin{eqnarray*}
c_{1\alpha\beta} & = & \frac{\partial t^i}{\partial\tau^1} 
\frac{\partial t^j}{\partial\tau^\alpha}
\frac{\partial t^k}{\partial\tau^\beta}
\left.c_{ijk}\right|_N \,, \\
& = & \left.c_{1jk}\right|_N
\frac{\partial t^j}{\partial\tau^\alpha}
\frac{\partial t^k}{\partial\tau^\beta}\,,\\
& = & \eta_{jk}
\frac{\partial t^j}{\partial\tau^\alpha}
\frac{\partial t^k}{\partial\tau^\beta}\,,\\
& = & \eta_{\alpha\beta}\,.
\end{eqnarray*}

\medskip

\noindent{\underline {Linearity of induced Euler vector field}}

\medskip

Let
\[
E=E^i \frac{\partial~}{\partial t^i}\,.
\]
Then, on using (\ref{decomposition})\,,
\[
\left.E\right|_\N = E^i 
\Big\{
A^\alpha_i \frac{\partial~}{\partial \tau^\alpha} +
n_i^{\tilde\alpha} \frac{\partial~}{\partial\nu^{\tilde\alpha}}
\Big\}\,.
\]
Thus, if $(\left.E\right|_\N)^\perp=0\,,$
\begin{eqnarray}
\left.E^i\right|_\N n_i^{~{\tilde\alpha}} & = & 0 \,,\label{eperpn} \\
E_\N^\alpha & = & \left.E^i\right|_\N \eta_{ij}\frac{\partial t^j}{\partial\tau^\beta}
\eta^{\alpha\beta}\,,\label{defen}
\end{eqnarray}
It follows from (\ref{eperpn}) that
\[
\left.E^i\right|_\N = \omega^\alpha \frac{\partial t^i}{\partial\tau^\alpha}
\]
for some function $\omega^\alpha(\tau)$ and with this (\ref{defen}) implies
$\omega^\alpha = E_\N^\alpha\,.$ Thus
\begin{equation}
E^i|_\N = E_\N^\alpha \frac{\partial t^i}{\partial\tau^\alpha}\,.
\label{eionn}
\end{equation}
To prove $E_\N$ is linear in $\tau$ its second derivatives will be calculated. From
(\ref{defen})
\[
\frac{\partial E_\N^\alpha}{\partial\tau^\sigma}  = \frac{\partial t^k}{\partial \tau^\sigma}
\left. \frac{\partial E^i}{\partial t^k}\right|_\N \eta_{ij}
\frac{\partial t^j}{\partial\tau^\beta} \eta^{\alpha\beta}+ \left.E^i\right|_\N \eta_{ij} \frac{\partial^2 t^j}{\partial\tau^\sigma\tau^\beta} \eta^{\alpha\beta}\,.
\]
Using (\ref{a4}) and (\ref{eperpn}), the second term vanishes.
Thus
\begin{equation}
\frac{\partial^2 E^\alpha_\N}{\partial\tau^\sigma\partial\tau^\nu}  = 
\eta_{ij} \frac{\partial^2 t^j}{\partial\tau^\beta\partial\tau^\nu}
\eta^{\alpha\beta} \frac{\partial t^k}{\partial\tau^\sigma}
\left. \frac{\partial E^i}{\partial t^k}\right|_\N+
\eta_{ij} \frac{\partial t^j}{\partial\tau^\beta}
\eta^{\alpha\beta} \frac{\partial^2 t^k}{\partial\tau^\sigma\partial\tau^\nu}
\left. \frac{\partial E^i}{\partial t^k}\right|_\N\,
\label{second}
\end{equation}
using the fact that $E^i$ is linear in $t\,.$ The first term in (\ref{second})
simplifies on using (\ref{a4}):
\[
{\rm first~term~} = \eta_{ij} \Omega_{\beta\nu}^{~~{\tilde\alpha}} n_{{\tilde\alpha}}^{~j}
\eta^{\alpha\beta} \frac{\partial t^k}{\partial\tau^\sigma}
\left. \frac{\partial E^i}{\partial t^k}\right|_\N\,.
\]
This simplified by first differentiating (\ref{eperpn}) and using (\ref{a6}), yielding
\[
{\rm first~term~} =
\eta^{\alpha\beta} \Omega_{\sigma\delta}^{~~{\tilde\alpha}} \Omega_{{\tilde\alpha}\beta\nu}
E_\N^\delta\,.
\]
The second term in (\ref{second}) may be written, using the explicit form $E^i=q^i_{~j}t^j+r^i$
and (\ref{a4}) as
\[
{\rm second~term~} = \eta_{ks}
\frac{\partial~}{\partial\tau^\beta}
\Big\{
\eta_{ij} q^i_{~r} \eta^{rs} t^j
\Big\}
\eta^{\alpha\beta} \Omega^{\tilde\alpha}_{~\sigma\nu} n^k_{~{\tilde\alpha}}\,.
\]
Using the explicit form $q^i_{~j}=(1-q_i) \delta_{ij}$ with $q_i+q_{m+1-i}=d\,,$
\[
\eta_{ij} q^i_{~r} \eta^{rs} = - q^s_{~j} + (2-d) \delta^s_j\,.
\]
Hence
\[
{\rm second~term~} = \Big\{
-\eta_{ks}\frac{\partial E^s}{\partial\tau^\beta} +
(2-d) \eta_{ks} \frac{\partial t^s}{\partial\tau^\beta}
\Big\} \eta^{\alpha\beta} \Omega^{\tilde\alpha}_{~\sigma\nu}\,.
\]
Repeating the earlier manipulations gives
\[
\frac{\partial^2 E_\N^\alpha}{\partial\tau^\sigma\partial\tau^\nu} =
\eta^{\alpha\beta} \eta^{{\tilde\alpha}{\tilde\beta}}
\big\{
\Omega^{\tilde\alpha}_{~\sigma\delta} \Omega^{\tilde\beta}_{~\beta\nu} -
\Omega^{\tilde\alpha}_{~\beta\delta} \Omega^{\tilde\beta}_{~\sigma\nu}
\big\}
E_\N^\delta
\]
and by virtue of the Gauss-Codazzi equation (\ref{a7}) this vanishes. Hence
$E_\N$ is linear in the $\tau$-variables.

\medskip

\noindent{\underline {Quasihomogeneity of induced prepotential}}

\medskip

The prepotential $F$ satisfies the quasihomogeneity condition
\[
{\cal L}_E F = d_F F + {\rm quadratic~terms~}\,.
\]
This is equivalent to the relation
\[
{\cal L}_E c_{ijk} = d_F\,c_{ijk}
\]
on structure functions. Expanding this gives
\[
E^r\frac{\partial c_{ijk}}{\partial t^r}
=d_F\,c_{ijk} - \frac{\partial E^r}{\partial t^i} c_{rjk} - {\rm cyclic}\,.
\]
Since the induced prepotential on $\N$ is only defined implicitly, the analogous relation
for the quasihomogeneity of the induced structure functions with respect to the induced
vector field will be found, the quasihomogeneity following by integration of this result.
The proof is straightforward:
\begin{eqnarray*}
E_\N c_{\alpha\beta\gamma} & = & E_\N^\sigma
\frac{\partial c_{\alpha\beta\gamma}}{\partial \tau^\sigma}\,,\\
& = & E_\N^\sigma \frac{\partial~}{\partial\tau^\sigma}
\Big\{
\frac{\partial t^i}{\partial \tau^\alpha}
\frac{\partial t^j}{\partial \tau^\beta}
\frac{\partial t^k}{\partial \tau^\gamma}
\left.c_{ijk}\right|_\N
\Big\}\,.
\end{eqnarray*}
But terms like
\[
E_\N^\sigma \frac{\partial^2 t^i}{\partial \tau^\sigma \partial \tau^\beta}
\frac{\partial t^j}{\partial \tau^\beta}
\frac{\partial t^k}{\partial \tau^\gamma}
\left.c_{ijk}\right|_\N = E_\N^\sigma \Omega_{\sigma\alpha}^{~~{\tilde\mu}}
\Xi_{\beta\gamma{\tilde\mu}}
\]
vanish since $\N$ is a natural submanifold. Thus
\begin{eqnarray*}
E_\N c_{\alpha\beta\gamma} & = &
\frac{\partial t^i}{\partial \tau^\alpha}
\frac{\partial t^j}{\partial \tau^\beta}
\frac{\partial t^k}{\partial \tau^\gamma}
E_\N^\sigma \frac{\partial t^p}{\partial \tau^\sigma}
\left. \frac{\partial c_{ijk}}{\partial t^p}\right|_\N\,.
\end{eqnarray*}
Using (\ref{eionn}) and the quasihomogeneity of $F$ gives
\[
E_\N(c_{\alpha\beta\gamma}) = d_F\,c_{\alpha\beta\gamma} -
\Big\{
\frac{\partial E^\sigma}{\partial\tau^\alpha} c_{\sigma\beta\gamma}+
E_\N^\sigma \Omega_{\alpha\sigma}^{~~{\tilde\alpha}} \Xi_{\beta\gamma{\tilde\alpha}}
\Big\}-{\rm cyclic}\,.
\]
Hence on a natural submanifold
\[
{\cal L}_{E_\N} c_{\alpha\beta\gamma} = d_F\,c_{\alpha\beta\gamma}\,,
\]
where ${\cal L}_{E_\N}$ is is Lie-derivative along $E_\N$ in the submanifold $\N\,.$
Integration then gives the quasihomogeneity of the induced prepotential. Note that
the total scaling dimension $d_F$ is unchanged.
This result is actually independent of the condition $\Xi=0\,,$ the terms involving $\Xi$
cancel. Thus on any submanifold where $(E_\N)^\perp=0$ the induced structure functions of
the not necessarily associative induced algebra are quasihomogeneous.

\medskip

\noindent This result may be formulated in terms of the vanishing of the induced
Dubrovin connection [D].

\bigskip

\subsection{The induced intersection form}

One important property of a Frobenius manifold is the existence of a second flat metric
defined by [D]
\begin{eqnarray*}
g^{ij}  & = & E(dt^i \circ dt^j) \,,\\
& = & c^{ij}_{~~k}E(dt^k)
\end{eqnarray*}
with the basic property that
\[
\frac{\partial g^{ij}}{\partial t^1} = \eta^{ij}\,.
\]
It follows from this that the pencil of metrics
\[
g_\lambda^{ij} = g^{ij} + \lambda \eta^{ij}
\]
if flat for all values of $\lambda\,.$ In this section it will be shown
(under the conditions of the above theorem) that the
restriction of this metric to the submanifold is given by the analogous
formulae. Since the above defines $g^{ij}$ rather than $g_{ij}\,,$
 a different approach is
required.

\medskip

Consider the tensor
\[
g^{ij} \frac{\partial~}{\partial t^i} \otimes \frac{\partial~}{\partial t^j}\,.
\]
Restricting this to $\N\,,$ and using (\ref{decomposition}) gives
\begin{eqnarray*}
\left.
g^{ij} \frac{\partial~}{\partial t^i} \otimes \frac{\partial~}{\partial t^j}
\right|_\N & = &
\left.g^{ij}\right|_\N \Big\{
A_i^{~{\alpha}} \frac{\partial ~}{\partial\tau^\alpha} +
n_i^{~{\tilde\alpha}} \frac{\partial ~}{\partial\tau^{\tilde\alpha}}
\Big\}
\otimes
\Big\{
A_j^{~{\beta}} \frac{\partial ~}{\partial\tau^\beta} +
n_j^{~{\tilde\beta}} \frac{\partial ~}{\partial\tau^{\tilde\beta}}
\Big\}\,, \\
& = &
\left.g^{ij}\right|_\N  A_i^{~{\alpha}}	A_j^{~{\beta}}
\frac{\partial~}{\partial\tau^\alpha} \otimes
\frac{\partial~}{\partial\tau^\beta} +
\left.g^{ij}\right|_\N
n_i^{~{\tilde\alpha}} n_j^{~{\tilde\beta}}
\frac{\partial ~}{\partial\tau^{\tilde\alpha}} \otimes
\frac{\partial ~}{\partial\tau^{\tilde\beta}} \\
& & \quad\quad\quad\quad\quad +2
\left.g^{ij}\right|_\N A_i^{~{\tilde\alpha}} n_j^{~{\tilde\beta}}
\frac{\partial~}{\partial\tau^\alpha} \otimes
\frac{\partial~}{\partial \nu^{\tilde\beta}}\,.
\end{eqnarray*}
Simply calculations show that, under the conditions of the above theorem,
\[
{\rm cross~term~} = 2 \eta^{\alpha\beta} \Xi_{\beta\sigma}^{~~{\tilde\beta}}
E^\sigma 
\frac{\partial~}{\partial\tau^\alpha} \otimes
\frac{\partial~}{\partial \nu^{\tilde\beta}}\,,
\]
and hence vanish. This gives an orthogonal decomposition and hence a metric
on $\N$ given by
\[
g^{\alpha\beta} = \left.g^{ij}\right|_\N  A_i^{~{\alpha}}	A_j^{~{\beta}}
\]
Similar manipulations give
\[
g^{\alpha\beta} = E_\N (d\tau^\alpha \star d \tau^\beta)\,.
\]
Thus the two ways to compute the induced intersection form, either
by the restriction of the intersection from on $\M$ to $\N\,,$ or by
calculating it using the induced Euler vector field on $\N\,$ agree.
Similarly
\[
\frac{\partial g^{\alpha\beta}}{\partial\tau^1} = \eta^{\alpha\beta}\,.
\]
One remaining question is to calculate the Weingarten operators for the
submanifold using this second metric.

\bigskip

\section{Frobenius submanifolds and the foldings of Coxeter graphs}

In this sections the above ideas will be applied to a class of Frobenius
manifolds constructed from a Coxeter group $W$ and in particular two dimensional
Frobenius submanifolds will be considered.

The full details of the construction of these Frobenius manifolds may be found in [D].
For these the Euler vector field takes the form
\[
E=\sum_{i=1}^m d_i t^i \frac{\partial~}{\partial t^i}\,,
\]
where the $d_i$ are the exponents of the Coxeter group, or equivalently, the degrees
of the basic $W$-invariant polynomials. These are given in Table [1] (Note the reverse
ordering, so $d_n=2\,,d_1=h\,.$) They satisfy the basic condition $d_i+d_{m+1-i}=h+2\,,$ where
$h$ is known as the Coxeter
number of the group.
The corresponding prepotential is polynomial,
and it has been conjectured that all such polynomial prepotentials arise from this
construction.

\begin{table}
\begin{tabular}{cccc|c}
~~~~~~~~~~~~~~~~~~~~~&&Coxeter Group&&Exponents
$d_n\,,\ldots\,,d_1=h$\\
\cline{2-5}
&&$A_n$&& $2\,,3\,,\ldots\,,n+1$ \\
&&$B_n$&& $2\,,4\,,6\,,\ldots\,, 2n$\\
&&$D_n$&& $2\,,4\,,6\,,\ldots\,, 2n-2\,,n$\\
&&$E_6$&& $2\,,5\,,6\,,8\,,9\,,12$\\
&&$E_7$&& $2\,,6\,,8\,,10\,,12\,,14\,,18$\\
&&$E_8$&& $2\,,8\,,12\,,14\,,18\,,20\,,24\,,30$\\
&&$F_4$&& $2\,,6\,,8\,,12$\\
&&$G_2$&& $2\,,6$\\
&&$H_3$&& $2\,,6\,,10$\\
&&$H_4$&& $2\,,12\,,20\,,30$\\
&&$I_2(m)$&& $2\,,m$
\end{tabular}
\caption{Degrees of the $W$-invariant polynomials. }
\end{table}

Using the parametrization (\ref{param}) together with the requirement that the induced metric
must be both flat and in flat coordinates implies that the two-dimensional submanifolds
are parametrized:

\begin{eqnarray*}
t^1 & = & \tau_1 - \frac{1}{2} \int \sum_{j=2}^{m-1} f_j^\prime(\tau_2)
f_{m+1-j}^\prime(\tau_2) d\tau_2\,,\\
t^j & = & f_j(\tau_2)\,, \quad j=2\,,\ldots\,, m-1\,, \\
t^m & = & \tau_2\,.
\end{eqnarray*}
If the condition $(\left.E\right|_\N)^\perp=0$ is now imposed one obtains simple equations
for the $f_i$ giving the parametrization
\begin{eqnarray*}
t^1 & = & \tau_1 - \frac{1}{4}
\left\{\sum_{j=2}^{m-1} k_j k_{m+1-j} d_j d_{m+1-j} \right\} \frac{1}{h} \tau_2^{h/2}\,,\\
t^j & = & k_j  \tau_2^{d_j/2}\,, \quad j=2\,,\ldots\,, m-1\,, \\
t^m & = & \tau_2\,
\end{eqnarray*}
(using the fact that $d_m=2$ for all Coxeter groups, remembering the reverse ordering of the
exponents) and the induced Euler vector field
\[
E_\N = h \tau^1 \frac{\partial~}{\partial\tau^1} + 2 \tau^2 \frac{\partial~}{\partial\tau^2}\,.
\]
By Proposition [2] this submanifold automatically is a Frobenius (but not necessarily natural)
submanifold and it is east to check that the induced prepotential is
\[
F_\N = \frac{1}{2} \tau_1^2 \tau_2 + p(k_i) \tau_2^{h+1}\,,
\]
where $p(k_i)$ is some function of the constants $k_i$ which define the submanifold.
This prepotential is polynomial and corresponds to the Coxeter group $I_2(h)\,.$ Thus
for any Coxeter group one has a family of two-dimensional Frobenius submanifold:
\[
\F_W \longrightarrow \F_{I_2(h)}\,.
\]
Natural Frobenius manifolds occur at special values of the constants $k_i\,.$

\medskip

\noindent{\bf Example}

Consider the Frobenius manifold
defined by the polynomial prepotential
\[
F_{H_3}=\frac{1}{2} t_1^2 t_3 + \frac{1}{2} t_1 t_2^2 + \frac{1}{60} t_2^3 t_3^2 +
\frac{1}{20} t_2^2 t_3^5 + \frac{1}{3960} t_3^{11}
\]
and Euler vector field
\[
E=10 t_1 \frac{\partial ~}{\partial t_1} + 6  t_2 \frac{\partial ~}{\partial t_2 }+
2 t_3 \frac{\partial ~}{\partial t_3}\,.
\]
Such a manifold is associated to the Coxeter group $H_3\,.$

By Proposition [2] any submanifold $\N$ defined by
\begin{eqnarray*}
t_1 & = & \tau_1 - \frac{9}{10} k_2^2 \tau_2^5 \,, \\
t_2 & = & k_2 \tau_2^3 \,, \\
t_3 & = & \tau_2
\end{eqnarray*}
is a Frobenius submanifold with respect to the induced structures. The condition for
the manifold to be a natural Frobenius
submanifold - normally a first order ordinary differential equation of degree three - reduces
to a cubic polynomial
\[
k_2(k_2-1)(27k_2+5)=0\,.
\]
Thus there are three natural Frobenius submanifolds of this form. This Frobenius submanifold
is also associated to
a Coxeter group, namely $I_2(10)\,.$ The relation between these two Coxeter groups may be seen
in terms of the folding of their Coxeter diagrams:

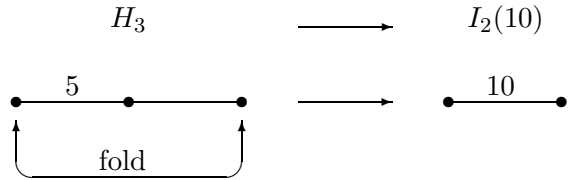
\setlength{\unitlength}{1.0cm}
\begin{picture}(11,3)
\put(4,1){\circle*{0.15}}
\put(4,1){\line(1,0){1.5}}
\put(5.5,1){\circle*{0.15}}
\put(5.5,1){\line(1,0){1.5}}
\put(7,1){\circle*{0.15}}
\put(7.75,1){\vector(1,0){1.25}}
\put(9.75,1){\circle*{0.15}}
\put(9.75,1){\line(1,0){1.5}}
\put(11.25,1){\circle*{0.15}}
\put(4.65,1.1){5}
\put(10.25,1.1){10}
\put(5.5,0.25){\oval(3,0.5)[b]}
\put(4,0.25){\vector(0,1){0.5}}
\put(7,0.25){\vector(0,1){0.5}}
\put(5.1,0.10){fold}
\put(7.75,2.0){\vector(1,0){1.25}}
\put(10.00,2.0){$I_2(10)$}
\put(5.25,2.00){$H_3$}
\end{picture}

\medskip

\noindent where such folding preserves the Coxeter number (in this case $10$) of the
groups involved. When $k_2=0$ the submanifold is just a plane, and for only this value
of $k_2$ does
\[
F_\N = \left.F\right|_\N\,.
\]
Similar results have been obtained by Zuber [Z] for natural Frobenius submanifolds obtained
by foldings of arbitrary Coxeter diagrams, but the only submanifolds that were considered
were hyperplanes. There are two other three-dimensional Coxeter groups,
namely $A_3$ and $B_3\,.$

\medskip

\noindent{\bf Example:} ${A_3\longrightarrow I_2(4)}$

The prepotential for the Frobenius manifold constructed from $A_3$ is
\[
F_{A_3} = \frac{1}{2} {t_1^2t_3}+\frac{1}{2} {t_1 t_2^2} + \frac{1}{4} {t_2^2t_3^2} +
\frac{1}{60} {t_3^5}\,.
\]
The two dimensional submanifold is given by
\begin{eqnarray*}
t_1 & = & \tau_1 - \frac{9}{16} k_2^2  \tau_2^2\,,\\
t_2 & = & k_2 \tau_2^{3/2}\,,\\
t_3 & = & \tau_2\,.
\end{eqnarray*}
The condition required for the submanifold to be a natural Frobenius submanifold reduce to
$k_2(32-27k_2^2)=0\,.$ 
Thus there are two natural Frobenius submanifolds given by $k_2=0\,,\pm\sqrt{32/27}\,,$
i.e. the plane $t_2=0$ and the cylinder over the semi-cubical parabola $27 t_2^2 = 32 t_3^3\,.$

\medskip

\noindent{\bf Example:} ${B_3\longrightarrow I_2(6)}$

The prepotential for the Frobenius manifold constructed from $B_3$ is
\[
F_{B_3} = \frac{1}{2} {t_1^2t_3}+\frac{1}{2} {t_1 t_2^2} + \frac{1}{6} {t_2^3t_3} +
\frac{1}{6} {t_2^2t_3^3} + \frac{1}{210} {t_3^7}\,.
\]
The two dimensional submanifold is given by
\begin{eqnarray*}
t_1 & = & \tau_1 - \frac{2}{3} k_2^2  \tau_2^3\,,\\
t_2 & = & k_2 \tau_2^2\,,\\
t_3 & = & \tau_2\,.
\end{eqnarray*}
The condition required for the submanifold to be a natural Frobenius submanifold reduce to
$k_2(2k_2-3)(-2k_2-1)=0\,.$ 
Thus there are three natural Frobenius submanifolds given by $k_2=0\,,-1/2\,,+3/2\,.$
\medskip

In these three examples the natural submanifolds are special from the point of view of
singularity theory, the submanifolds are cylinders over the caustics of $A_3\,,B_3$ and
$H_3\,.$ This observation does not generalize directly, for example the cylinder over the
caustic of $A_4$ is not a flat submanifold, and so cannot be a Frobenius submanifold.
However the induced multiplication is associative and quasihomogeneous
(since $(E_\N)^\perp=0$). These properties are best understood in terms of weak Frobenius and
F-manifolds [H,HM].

\medskip

\noindent{\bf Example:} $F_4\longrightarrow I_2(12)$

As a higher dimensional example, consider the embeddings of $I_2(12)$ in $F_4\,.$ The
prepotential for the Frobenius manifold constructed from $F_4$ is
\[
F_{F_4}=\frac{1}{2} {t_1^2 t_4} +t_1 t_2 t_3+ \frac{1}{6} {t_2^3 t_4} +
\frac{1}{12} t_3^4 t_4 + 
\frac{1}{6} {t_2 t_3^2 t_4^3} + \frac{1}{60} {t_2^2 t_4^5} + \frac{1}{252} {t_3^2 t_4^7} + 
\frac{1}{185328} {t_4^{13}}\,.
\]
The two dimensional submanifold is given by
\begin{eqnarray*}
t_1 & = & \tau_1 - 2 k_2 k_3 \tau_2^6\,,\\
t_2 & = & k_2 \tau_2^4\,,\\
t_3 & = & k_3 \tau_2^3\,,\\
t_4 & = & \tau_2\,.
\end{eqnarray*}
The conditions required for the submanifold to be a natural Frobenius submanifold are
\begin{eqnarray*}
k_2+12k_2^2+5 k_3^2 - 36 k_2 k_3^2 & = & 0 \,,\\
k_3(1+36k_2-144k_2^2+36k_3^2) & = & 0\,.
\end{eqnarray*}
These algebraic equations are easily solved giving six two-dimensional
natural Frobenius submanifolds (ignoring one complex solution):

\[
(k_2,k_3) = \left\{
	\begin{array}{ll}
	(0,0)\,, & (-1/12,0)\,,\\
        (-1/36,+1/18)\,, & (+5/12,+1/2)\,,\\   
        (-1/36,-1/18)\,,  & (+5/12,-1/2)\,.
        \end{array}
\right.
\]
Further examples may easily be constructed using the known formulae
for prepotentials constructed from Coxeter groups [Z].

\section{The quantum cohomology of $\mathbb{CP}^2$}

The quantum cohomology of
$\mathbb{CP}^2$ is given in terms of the prepotential 
\[
F=\frac{1}{2}t_1^2 t_3 + \frac{1}{2} t_1 t_2^2 +
\sum_{n=1}^\infty
\frac{N_n^{(0)} t_3^{3n-1} e^{nt_2}}{(3n-1)!}\,,
\]
with
\[
E=t_1\frac{\partial ~}{\partial t_1} + 3 \frac{\partial ~}{\partial t_2}-
t_3 \frac{\partial ~}{\partial t_3}\,,
\]
where $N_n^{(0)}$ is the number of rational curves of degree $n$ through $3n-1$ generic
points. The equations of associativity imply the recursion relation
\[
N_n^{(0)} = \sum_{i+j=n}
\Bigg[ \pmatrix{ 3n-4 \cr 3i-2} i^2 j^2 - i^3 j \pmatrix{ 3n-4 \cr 3i-1}\Bigg]
N_i^{(0)} N_j^{(0)}
\]
first derived by Kontsevich and Manin. With the initial condition
$N_1^{(0)}=1$ this determines all the $N_n^{(0)}\,.$ Following the
derivation in [D], this prepotential may be written as
\[
F=\frac{1}{2}t_1^2 t_3 + \frac{1}{2} t_1 t_2^2 + t_3^{-1} \phi(x)
\,,
\]
where $x=t_2 + 3 \log t_3\,.$ The equations of associativity then reduce to the
third order ordinary differential equation
\begin{equation}
9\phi^{'''}-18 \phi^{''}+11\phi^{'} - 2 \phi = \phi^{''} \phi^{'''} - \frac{2}{3}
\phi^{'} \phi^{''} + \frac{1}{3} {\phi^{''}}^2\,,
\label{ode}
\end{equation}
and with the ansatz
\begin{equation}
\phi(x) = \sum_{n=1}^\infty \frac{N_n^{(0)} }{ (3n-1)!} e^{nx}
\label{series}
\end{equation}
one obtains the above recursion relation.

By Proposition [2] any suitable two dimensional submanifold is a Frobenius manifold, but
a particularly interesting submanifold is given by $x=x_0\,,$ where $x_0$ is a constant.
On such a submanifold $(\left.E\right|_\N)^\perp=0\,.$ In terms of the parametrization
\begin{eqnarray*}
t_1 & = & \tau_1 + \frac{9}{2} \tau_2^{-1}\,,\\
t_2 & = & x_0 - 3 \log \tau_2 \,, \\
t_3 & = & \tau_2
\end{eqnarray*}
one obtains a Frobenius manifold on $\N$ defined by
\begin{eqnarray*}
F_\N & = & \frac{1}{2} \tau_1^2 \tau_2 - \Bigg[ \frac{ 81-8\phi(x_0) + 20 \phi^{'}(x_0)}{8}
\Bigg] \, \tau_2^{-1}\,,\\
E_\N & = &
\tau_1 \frac{\partial ~}{\partial \tau_1} - \tau_2 \frac{\partial ~}{\partial \tau_2}\,.
\end{eqnarray*}
The obstruction to this being a natural Frobenius submanifold is
\begin{equation}
27 + 2 \phi^{'}(x_0) - 3 \phi^{''} (x_0) = 0\,.
\label{obstructionCP2}
\end{equation}
It is not immediately obvious that a natural submanifold exists.

\bigskip

\begin{lemma} There exists a natural Frobenius submanifold, given by the condition
$x=x_0\,,$ where $x_0$ is the radius of convergence of the series (\ref{series}).
\end{lemma}

\bigskip

\noindent{\bf Proof~} It was shown in [FI] that the series (\ref{series}) has a finite
radius of convergence $x_0\,.$ Moreover it was shown that
$\phi\,,\phi^{'}\,,\phi^{''}\,,\phi^{'''}$ are all positive with
$\phi<\phi^{'}<\phi^{''}<\phi^{'''}$ for real $x<x_0\,,$ and that
$\phi\,,\phi^{'}$ and $\phi^{''}$ remains finite at $x_0$ with $\phi^{'''}$ blowing up.
Using these results, in the vicinity of $x_0\,,\phi$ takes the form
\[
\phi=\phi_0 + \phi_1(x_0-x) + \phi_2 \frac{(x_0-x)^2}{2} + \lambda (x_0-x)^{\alpha+2} +
\ldots\,,
\]
and substituting this into the differential equation (\ref{ode}) and equating coefficients
yields $\alpha=1/2\,,\lambda$ and the relation (\ref{obstructionCP2})\,. Hence a
natural Frobenius submanifolds exists.

\bigskip

\bigskip

\section{The quantum cohomology of $\mathbb{CP}^1 \times \mathbb{CP}^1$}

As explained elsewhere [FI], the quantum cohomology of
$\mathbb{CP}^1\times\mathbb{CP}^1$ is given in terms of the prepotential
\[
F=\frac{1}{2} t_1^2 t_4 + t_1 t_2 t_3 +
\sum_{\scriptstyle a\,,b \geq 0 \atop \scriptstyle a+b \geq 1}
\frac{N_{ab}}{[2(a+b)-1]!} t_4^{2(a+b)-1} e^{at_2+b t_3}\,,
\]
and Euler vector field
\[
E=\frac{1}{2} t_1 \frac{\partial ~}{\partial t_1}+\frac{\partial ~}{\partial t_2}+
\frac{\partial ~}{\partial t_3}- \frac{1}{2} t_4 \frac{\partial ~}{\partial t_4}\,.
\]
The coefficients $N_{ab}$ are the number of rational curves on a smooth quadric (such
quadrics being isomorphic to
$\mathbb{CP}^1 \times \mathbb{CP}^1$)
with bidegree $(a,b)$ though $2(a+b)-1$ points. These are determined by the initial
conditions $N_{01}=1\,,N_{ab}=N_{ba}$ and the recursion relations
\begin{eqnarray*}
2abN_{ab}&=&\sum N_{a_1b_1}N_{a_2b_2}a_1^2b_2^2(a_1b_2-a_2b_1)
\pmatrix{ 2(a+b)-2 \cr 2(a_1+b_1)-1 }\,,\\
aN_{ab}&=&\sum N_{a_1b_1}N_{a_2b_2}a_1(a_1^2b_2^2-a_2^2b_1^2)
\pmatrix{ 2(a+b)-3 \cr 2(a_1+b_1)-1 }\,,\\
0&=&\sum N_{a_1b_1}N_{a_2b_2}a_1^2
[(a_2+b_2-1)(b_1a_2+b_2a_1)-a_2b_2(2(a_1+b_1)-1)]
\pmatrix{ 2(a+b)-3 \cr 2(a_1+b_1)-1 }\,,\\
N_{ab}&=&\sum N_{a_1b_1}N_{a_2b_2}(a_1b_2+a_2b_1)b_2
\Bigg[
a_1\pmatrix{2(a+b)-4 \cr 2(a_1+b_1)-2 }-
a_2\pmatrix{2(a+b)-4 \cr 2(a_1+b_1)-3 }
\Bigg]\,,
\end{eqnarray*}
where the sums are over $a_1\,,a_2\,,b_1\,,b_2\geq 0\,,a_1+a_2=a\,,b_1+b_2=b\,.$

The symmetry $t_2 \longleftrightarrow t_3$ in these formulae suggest that one should
consider the codimension one submanifold defined by the parametrization
\begin{eqnarray*}
t_1 & = & \tau_1 \,, \\
t_2 & = & \frac{1}{\sqrt{2}} \tau_2 \,, \\
t_3 & = & \frac{1}{\sqrt{2}} \tau_2 \,, \\
t_4 & = & \tau_3\,,
\end{eqnarray*}
where the factor $\sqrt{2}$ ensures that the induced metric takes the canonical antidiagonal
form. This submanifold also satisfies the condition $(\left.E\right|_N)^\perp=0$ so
\[
E_\N = \frac{1}{2} \tau_1 \frac{\partial ~}{\partial\tau_1} +
\sqrt{2} \frac{\partial ~}{\partial\tau_2} - \frac{1}{2}\tau_3\frac{\partial ~}{\partial\tau_3}\,.
\]
The calculation of the induced multiplicaion on $\N$ is particularly simple, due to the fact
that $\N$ is just a hyperplane. The induced structure $\F_\N$ is a natural Frobenius
submanifold, the obstructions all take the form
\[
\Xi = \sum(a-b) S(a,b)
\]
with $S(a,b)=S(b,a)$ and hence vanish. The induced prepotential is given by
\begin{eqnarray*}
F_\N&=& \left.F\right|_\N \,, \\
& = & \frac{1}{2} \tau_1^2 \tau_3 + \frac{1}{2} \tau_1 \tau_2^2 + \tau_3^{-1}
\sum_{n=1}^\infty
\frac{ \Big[ \sum_{r=0}^n N_{n-r,r} \Big] }{(2n-1)!} \tau_3^{2n} e^{n\tau_2/\sqrt{2}}\,.
\end{eqnarray*}

While this construction guarantees that $\F_\N$ is a Frobenius manifold it is interesting
to calculate the relations required to ensure that the prepotential
\[
F = \frac{1}{2} \tau_1^2 \tau_3 + \frac{1}{2} \tau_1 \tau_2^2 + \tau_3^{-1}
\sum_{n=1}^\infty
\frac{ N_n }{(2n-1)!} \tau_3^{2n} e^{n\tau_2/\sqrt{2}}
\]
defines a Frobenius manifold. The calculations are identical, apart from
different numerical coefficients, to the calculation of the quantum cohomology of
$\mathbb{P}^2$ so the details will not be repeated. It turns out that the coefficients
$N_n$ must satisfy the recursion relation
\[
N_n=\frac{1}{2} (2n-4)! \sum_{\scriptstyle k\geq 1\,,l\geq 1 \atop \scriptstyle k+l=n}
\frac{kl[kl(n+1) - (l^2+k^2))]}{(2k-1)!(2l-1)!} N_k N_l
\]
with initial condition $N_2=2\,.$ Thus the numbers $N_n=\sum_{r=0}^n N_{n-r,r}$ must satisfy
the above recursion relation. This may be easily verified for small values of $n$, but the
fact that $\N$ is a {\sl natural} Frobenius submanifold makes the result automatic. Presumably
one may also derive this result directly from the recursion relations
Obviously the numbers $N_n$ contain less information the the original $N_{ab}\,,$ the Frobenius
submanifold only determining their sum, not the individual numbers.

\begin{table}
\begin{tabular}{cccc|c}
~~~~~~~~~~~~~~~~~~~~~&&$n$&&$N_n=\sum_{r=0}^n N_{n-r,r}$\\
\cline{2-5}
&&1&&2 \\
&&2&&1 \\
&&3&&2 \\
&&4&&14 \\
&&5&&194 \\
&&6&&4792 \\
&&7&&182770 \\
&&8&&10078480 \\
&&9&&758120642 \\
&&10&&74795167616 \\
&&11&&937456239394 \\
&&12&&1456089241205248
\end{tabular}
\caption{The numbers $N_n$ for $1\leq n \leq 12$} 
\end{table}

One may also, mirroring the construction in the last section, obtain a Frobenius
submanifold of $\F_\N$ on the submanifold of $\N$ defined by the condition
\[
\frac{1}{\sqrt{2}} \tau_2 + 2 \log \tau_3 = {\rm constant}\,.
\]
Thus one obtains a nested sequence of Frobenius manifolds.

Underlying this construction is the symmetry $t_2\longleftrightarrow t_3\,.$ The
origin of this symmetry comes from the fact that the Frobenius manifold is a
tensor product of two 2-dimensional Frobenius manifolds [K],
\begin{equation}
\F_{\mathbb{CP}^1\times\mathbb{CP}^1} \cong \F_{\mathbb{CP}^1} \otimes \F_{\mathbb{CP}^1}\,,
\label{product}
\end{equation}
where $\F_{\mathbb{CP}^1}$ is given by
\begin{eqnarray*}
F_{\mathbb{CP}^1} & = & \frac{1}{2} t_1 t_2^2 + e^{t_2}\,,\\
E_{\mathbb{CP}^1} & = & \frac{1}{2} t_1
\frac{\partial ~}{\partial t_1} + 2 \frac{\partial ~}{\partial t_2}\,.
\end{eqnarray*}
The Euler vector field for the
product (\ref{product})\,, constructed from $E_{\mathbb{CP}^1}\,,$
is
\[
E_{\mathbb{CP}^1\times\mathbb{CP}^1} =
t^{11} \frac{\partial ~}{\partial t^{11}} + 2
t^{12} \frac{\partial ~}{\partial t^{12}} + 2
t^{21} \frac{\partial ~}{\partial t^{21}} -
t^{22} \frac{\partial ~}{\partial t^{22}}
\]
and this, by construction, automatically has the required
symmetry $t^{12}\longleftrightarrow t^{21}\,.$
Thus the natural Frobenius submanifold may be formulated in terms of a quotient of this
product by this symmetry:
\[
\F_\N \cong \frac{
\F_{\mathbb{CP}^1} \times \F_{\mathbb{CP}^1} }{ t^{12}\longleftrightarrow t^{21}}\,.
\]
More generally one may obtain new Frobenius manifolds by squaring a Frobenius manifold and
taking such a quotient
\[
\F_\N \cong \frac{\F_\M \otimes \F_\M}{\sim}\,.
\]

\bigskip

\noindent{\bf Example} Another example of this kind is given in terms
of the Frobenius manifold
\[
F_{A_2}=\frac{1}{2} t_1^2 t_2 + t_2^4\,
\]
which is constructed from the Coxeter group $A_2 \cong I_2(3)\,.$ The
product of two such manifolds is again a Frobenius manifold associated
to the Coxeter group $D_4\,:$
\[
\F_{D_4} \cong \F_{A_2} \otimes \F_{A_2}\,.
\]
By construction this automatically has the symmetry
$t^{12} \longleftrightarrow t^{21}$ so, as in the case of
$\F_{\mathbb{CP}^1} \otimes \F_{\mathbb{CP}^1}$ one has a natural Frobenius submanifold
defined on the hyperplane $t^{12}=t^{21}\,.$ This Frobenius submanifold is again associated
to a Coxeter group, namely $B_3\,,:$
\[
\F_{B_3} \cong \frac{\F_{A_2} \otimes \F_{A_2}}{\sim}\,.
\]
Repeating the construction outlined in section 1, one obtains natural Frobenius submanifolds
inside $\F_{B_3}\,,$ this time associated to the Coxeter group $I_2(6)\,.$ Thus one obtains
a nested sequence of natural Frobenius manifolds
\[
\F_{I_2(6)} \subset \F_{B_3} \subset \F_{A_2} \otimes \F_{A_2} \cong \F_{D_4}\,,
\]
This sequence may be understood in terms of foldings of Coxeter diagrams:

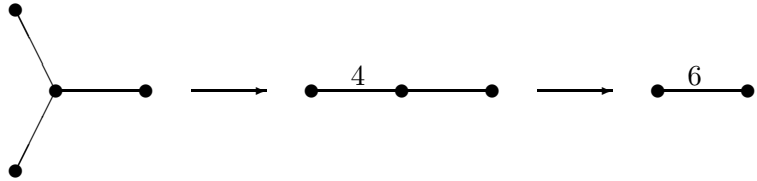
\setlength{\unitlength}{0.8cm}
\begin{picture}(15,3)
\put(7,1){\circle*{0.2}}
\put(7,1){\line(1,0){1.5}}
\put(8.5,1){\circle*{0.2}}
\put(8.5,1){\line(1,0){1.5}}
\put(10,1){\circle*{0.2}}
\put(10.75,1){\vector(1,0){1.25}}
\put(12.75,1){\circle*{0.2}}
\put(12.75,1){\line(1,0){1.5}}
\put(14.25,1){\circle*{0.2}}
\put(7.65,1.1){4}
\put(13.25,1.1){6}

\put(5,1){\vector(1,0){1.25}}
\put(4.25,1){\circle*{0.2}}
\put(2.75,1){\circle*{0.2}}
\put(2.75,1){\line(1,0){1.5}}
\put(2.75,1){\line(-1,2){0.67}}
\put(2.75,1){\line(-1,-2){0.67}}

\put(2.08,2.34){\circle*{0.2}}
\put(2.08,-0.34){\circle*{0.2}}

\end{picture}

\bigskip

\noindent One may also embed the trivial 1-dimensional Frobenius manifold
given by $F=t_1^3/6\,$ in $F_{I_2(6)}\,,$ giving complete nested sequence
of Frobenius submanifolds.

\section{Conclusion}
The results of this paper have been derived using flat-coordinates only. One
important avenue for future research is to rederive them using canonical
coordinates. Such an approach will involve the classical differential geometric problem
of properties of flat submanifolds of Ergoff metrics which are themselves Ergoff.
One basic object that is best studies using canonical coordinates is the isomonodromic
$\tau$-function, denoted by $\tau_I\,.$ One obvious question is how the $\tau_I$-function
of  a (natural)-submanifold is related to that of its parent Frobenius manifold. As the
following discussion will show, the relation, whatever it is, is not straightforward.

One way to study certain properties of the $\tau_I$-function without having to use canonical
coordinates is to use the following result
\begin{equation}
\tau_I = J^\frac{1}{24} \, e^G
\label{DZtauI}
\end{equation}
recently proved in [DZ] for semi-simple Frobenius manifolds. Here $J$ is the Jacobian of the
transformation from canonical to flat coordinates, and $G$ is the solution to Getzler's
equations for genus-one Gromov-Witten invariants. Consider the Frobenius manifolds
constructed from the Coxeter groups $A_3\,,B_3$ and $H_3\,.$ The corresponding $G$-functions
are
\begin{eqnarray*}
G_{A_3} & = & 0 \,, \\
G_{B_3} & = & -\frac{1}{48} \log [2t_2-3t_3^2]\,, \\
G_{H_3} & = & -\frac{1}{20} \log [t_2 - t_3^3]\,.
\end{eqnarray*}
In the later two cases $G$ has a logarithmic singularity on {\sl one} of the corresponding
natural Frobenius submanifolds. In all cases these natural submanifolds lie in the nilpotent
locus, so from (\ref{DZtauI}) the $\tau_I$-function is singular on all three of the natural
Frobenius submanifolds. This property is also present in the Frobenius manifold for the quantum
cohomology of $\mathbb{CP}^2\,.$ The derivative of the $G$-function is given by
\[
G^\prime = \frac{\phi^{\prime\prime\prime}-27}{8(27+2\phi^\prime-3\phi^{\prime\prime})}\,,
\]
where $\phi$ is given by (\ref{ode}), and using the series expansion (\ref{series}) one may
integrate this equation and show that $G$ also has a logarithmic singularity on the natural
Frobenius submanifold. This submanifold does not lie in the nilpotent locus, and it follows
from (\ref{DZtauI}) that $\tau_I$ is also singular on the submanifold.

For the quantum cohomology of $\mathbb{CP}^1 \times\mathbb{CP}^1$ it is unclear what the
singularity structure of the $G$-functions is since its governing equations are somewhat
more complicated, but even if the $G$-function does restrict to the submanifolds, it does
not restrict to the $G$-function of the submanifold. This is easily seen by calculating the
scaling constant $\gamma$ defined by ${\cal L}_E G = \gamma\,.$ The scaling constant of the
$\left.G\right|_\N$ does not equal the scaling constant of $G_\N\,.$

It has been shown that the contracted genus-zero Gromov-Witten invariants
$\sum_{a+b=n} N_{ab}^{(0)}$ satisfy a simple recursion relation which may be understood
as coming from a natural codimension one Frobenius submanifold. This raises the question of how
higher-genus contracted Gromov-Witten invariants $\sum_{a+b=n} N_{ab}^{(g)}$ are related, if
at all, to this submanifold. It would also be of interest both to have a direct proof of the
genus zero result by contracting the full recursion relations for the $N_{ab}^{(0)}$,
and to have a
direct algebraic-geometric proof of why this submanifold \lq counts\rq~these contracted
sums.

In summary, the results suggest the following problems:

\begin{itemize}

\item[$\bullet$] How can one reformulate these results in terms of canonical coordinates?

\item[$\bullet$] How is the singularity structure of the $G$-function related to the
existence of natural Frobenius submanifolds?

\item[$\bullet$] If $\N\subset\M$ is a natural Frobenius submanifold, what are the
relationships
\begin{eqnarray*}
(\tau_I)_\N & \longleftrightarrow & (\tau_I)_\M \,, \\
G_\N & \longleftrightarrow & G_\M ?
\end{eqnarray*}
These are clearly related by (\ref{DZtauI}). For the KP hierarchy there are some
interesting results [AvM] on the Birkhoff strata of the Grassmannian based on the zeros of the
$\tau$-function. It would be interesting to study the dispersionless counterparts of such
systems.

\end{itemize}

\noindent Finally, it should be possible to study degenerate Frobenius manifolds introduced
in [S2] in this framework, by embedding them in higher-dimensional, non-degenerate Frobenius
manifolds [Ko].

\section*{Acknowledgments}

I would like to thank Claus Hertling for his comments on this work, and in particular for
pointing out the relation between caustics and Frobenius submanifolds.

\appendix

\section*{Appendix}

\renewcommand{\theequation}
{A.\arabic{equation}}
\setcounter{equation}{0}

Since a Frobenius manifold is flat, any Frobenius submanifold must also be flat, and hence one has
to consider the possible embedding of one flat space in another. The following results are
entirely standard (see for example [E]) and are just a specialization of the
general Gauss-Codazzi equations for the embedding of an arbitrary manifold into another.

From (\ref{basisN}) the induced metric on $\N$ is
\begin{equation}
\eta_{\alpha\beta}= \eta_{ij} \frac{\partial t^i}{\partial\tau^\alpha}
\frac{\partial t^j}{\partial\tau^\beta}\,.
\label{a1}
\end{equation}
It will be assumed that the $\tau^\alpha$-coordinates are flat coordinates, i.e. the components
of $\eta_{\alpha\beta}$ are constant. Let $n^j_{~{\tilde\alpha}}$ be a field of normal vectors
to $\N\,,$ so
\begin{equation}
\eta_{ij} \frac{\partial t^i}{\partial\tau^\alpha} n^j_{~{\tilde\alpha}} = 0 \,,
\label{a2}
\end{equation}
normalized so
\begin{equation}
\eta_{ij} n^i_{~{\tilde\alpha}} n^j_{~{\tilde\beta}}=\eta_{{\tilde\alpha}{\tilde\beta}}
\label{a3}
\end{equation}
where $\eta_{{\tilde\alpha}{\tilde\beta}}$ are constant with
$\eta_{{\tilde\alpha}{\tilde\beta}}=\epsilon({\tilde\alpha})
\delta_{{\tilde\alpha}{\tilde\beta}}$ with $\epsilon({\tilde\alpha})=\pm 1\,.$

Differentiating (\ref{a1}) implies
\[
\eta_{ij} \frac{\partial^2 t^i}{\partial\tau^\alpha\partial\tau^\beta} \,
\frac{\partial t^j}{\partial\tau^\beta} = 0\,,
\]
and hence there exist functions $\Omega_{\alpha\beta}^{~~~{\tilde\alpha}}$ such that
\begin{equation}
\frac{\partial^2 t^i}{\partial\tau^\alpha\partial\tau^\beta}=
\Omega_{\alpha\beta}^{~~~{\tilde\alpha}} n_{\tilde\alpha}^{~ i}\,.
\label{a4}
\end{equation}
Differentiating (\ref{a2}) implies, on using (\ref{a4})
\begin{equation}
\Omega_{{\tilde\alpha}\alpha\beta} = - \eta_{ij}
\frac{\partial t^i}{\partial\tau^\alpha}\,
\frac{\partial n_{{\tilde\alpha}}^{~ j}}{\partial\tau^\beta}\,.
\label{a5}
\end{equation}
Differentiating (\ref{a3}) implies that
\begin{equation}
\frac{\partial n^j_{~{\tilde\beta}}}{\partial\tau^\alpha} = - \Omega_{{\tilde\beta}\alpha\sigma}
\frac{\partial t^j}{\partial\tau^\nu} \eta^{\sigma\mu}\,.
\label{a6}
\end{equation}
Note in particular that the torsion tensors are zero. The immediate consequence of this is that
the normal bundle of $\N$ is flat, i.e. $d {\vec n}_\alpha \in T\N\,.$

The Gauss-Codazzi equations, the integrability conditions for the above structures, reduce to the
three equations
\begin{equation}
\eta^{{\tilde\alpha}{\tilde\beta}}\big[
\Omega_{{\tilde\alpha}\alpha\beta} \Omega_{{\tilde\beta}\gamma\delta} -
\Omega_{{\tilde\alpha}\alpha\delta} \Omega_{{\tilde\beta}\gamma\beta}\big] = 0\,,
\label{a7}
\end{equation}
and
\begin{equation}
\eta^{\mu\nu}\big[
\Omega_{{\tilde\alpha}\mu\alpha} \Omega_{{\tilde\beta}\nu\beta}-
\Omega_{{\tilde\alpha}\mu\beta} \Omega_{{\tilde\beta}\nu\alpha}\big]=0\,,
\label{a8}
\end{equation}
and
\begin{equation}
\frac{\partial \Omega_{{\tilde\alpha}\alpha\mu}  }{\partial\tau^\nu} -
\frac{\partial \Omega_{{\tilde\alpha}\alpha\nu}  }{\partial\tau^\mu} = 0\,.
\label{a9}
\end{equation}

\section*{Bibliography}

\begin{itemize}

\item[{[AvM]}] {Adler, M. and van Moerbeke, P., {\sl Birkhoff Strata, B\"acklund transformations
and regularization of isospectral operators} Adv. Math. {\bf 108} (1994) 140-204.}

\item[{[D]}] {Dubrovin, B., {\sl Geometry of 2D topological field theories} in {\sl Integrable
Systems and Quantum Groups}, Editors: M. Francaviglia and S. Greco. Springer lecture
notes in mathematics, {\bf 1620}, 120-348.}

\item[{[DZ]}] {Dubrovin, B. and Zhang, Y.,
{\sl Bihamiltonian hierarchies in the 2D Topological
Field Theory at One-Loop Approximation}, C.M.P. {\bf 198} (1998) 311-361,
{\sl Frobenius Manifolds and Virasoro Constraints}, math/9808048.}

\item[{[E]}] {Eisenhart, L.P. {\sl Riemannian geometry}, Princeton Univ. Press (1949).}

\item[{[FI]}] {Di Francesco, P. and Itzykson, C., {\sl Quantum intersection rings}
hep-th/9412175.}

\item[{[H]}] {Herling, C., {\sl Multiplication on the tangent bundle}, math/9910116.}

\item[{[HM]}] {Hertling, C. and Manin, Yu., {Weak Frobenius manifolds},
Int. Math. Res. Notices {\bf 6} (1999) 277-286.}

\item[{[K]}] {Kaufmann, R.M., {\sl The tensor product in the theory of Frobenius manifolds},
Int. J. Math. {\bf 10:2} (1999) 159-206.}

\item[{[Ko]}] {Kodama, Y., {\sl Dispersionless integrable systems and their solutions}, to appear
in the proceedings of the conference
{\sl Integrability: the Seiberg-Witten and Whitham equations}
held at the ICMS Edinburgh in 1998.}

\item[{[S1]}] {Strachan, I.A.B., {\sl Jordan manifolds and dispersionless KdV equations},
in preparation.}

\item[{[S2]}] {Strachan, I.A.B., {\sl Degenerate Frobenius manifolds and the bi-Hamiltonian
structure of rational Lax equations}, J. Math. Phys. {\bf 40:10} (1999) 5058-5079.}

\item[{[Z]}] {Zuber, J.-B., {\sl On Dubrovin topological field theories}, Mod. Phys. Lett.
{\bf A9} (1994) 749-760.}

\end{itemize}

\end{document}